\documentclass[12pt]{article}

\usepackage{graphicx}
\usepackage{pst-poly}     
\usepackage{pst-all}      
\usepackage{amsfonts,amssymb,amsmath,latexsym,bm}

\newtheorem{thm}{Theorem}[section]
\newtheorem{Def}[thm]{Definition}

\newtheorem{pps}[thm]{Proposition}
\newtheorem{cor}[thm]{Corollary}
\newtheorem{lem}[thm]{Lemma}

\newenvironment{pf}[1][Proof]{\noindent\textbf{#1.} }{\hfill\rule{1mm}{2mm}}

\makeatletter \@addtoreset{equation}{section} \makeatother

\begin{document}

\title{Pebbling on Jahangir graphs
  }
\author
{Zheng-Jiang Xia\footnote{Corresponding
author:
xiazj@aufe.edu.cn}, Zhen-Mu Hong, Fu-Yuan Chen\\ \\
{\small     School of Finance,} \\
{\small             Anhui University of Finance and Economics.}   \\
\small Bengbu, Anhui, 233030, China\\
\small Email:xzj@ustc.edu.cn, zmhong@mail.ustc.edu.cn, 757552169@qq.com\\
}
\date{}
\maketitle

\begin{quotation}
\noindent\textbf{Abstract:}The  pebbling  number  of  a  graph  $G$, $f(G)$,  is  the  least  $p$  such  that,  however  $p$  pebbles  are  placed  on  the  vertices  of  $G$,  we  can  move  a  pebble  to  any  vertex by  a  sequence  of  moves, each  move  taking  two  pebbles  off one  vertex  and  placing one  on  an  adjacent  vertex. In this paper, we will show the pebbling number of Jahangir graphs $J_{n,m}$ with $n$ even, $m\geq8$.

\end{quotation}



\section{Introduction}

Pebbling in graphs was first introduced by Chung\cite{c89}. We only consider simple connected graphs. For a given graph $G$, the pebbling distribution $D$ of $G$ is a projection from $V(G)$ to $N$, $D: V(G)\rightarrow N$, where, $D(v)$ is the number of pebbles on $v$, the total number of pebbles on a subset $A$ of $V$ is given by $|D(A)|=\Sigma_{v\in A}D(v)$, $|D|=|D(V)|$ is the size of $D$.

A \emph{pebbling move} consists of the removal of two pebbles from a vertex and the placement of one pebble on an adjacent vertex.
Let $D$ and $D'$ be two pebbling distribution of $G$, we say that $D$ contains $D'$ if $D(v)\geq D'(v)$ for all $v\in V(G)$, we say that $D'$ is reachable from $D$ if there is some sequence (probably empty) of pebbling moves start from $D$
and resulting in a distribution that contains $D'$.
For a graph $G$, and a vertex $v$, we call $v$ a \emph{root} if the goal is to place pebbles on $v$; If $t$ pebbles can be moved to $v$ from $D$ by a sequence of pebbling moves, then we say that $D$ is $t$-fold \emph{$v$-solvable}, and $v$ is $t$-\emph{reachable} from $D$. If $D$ is $t$-fold $v$-solvable for every vertex $v$, we say that $D$ is $t$-\emph{solvable}.

Computing the pebbling number is difficult in general. The problem of deciding if a given distribution on a graph can reach a particular vertex
was shown in \cite{mc06} to be NP-complete, even for planar graphs\cite{lcd14}.  The problem of deciding whether a graph $G$ has pebbling number at most $k$ was shown in \cite{mc06} to be $\Pi^P_2$-complete.


%
%
%
%
%
%
Given a root vertex $v$ of a tree $T$, then we can view $T$ be a directed graph $\vec{T}_v$ with each edge directed to $v$, A path
partition is a set of nonoverlapping directed paths the union of which is  $\vec{T}_v$. A path partition
is said to \emph{majorize }another if the nonincreasing sequence of the path size majorizes
that of the other (that is $(a_1,a_2,\ldots,a_r)>(b_1,b_2,\ldots,b_t)$ if and only if $a_i > b_i$  where $i= \min\{j:a_j\neq b_j \}$). A path partition of a tree $T$ is said to be \emph{maximum} if it majorizes all other path partitions.

\begin{thm}{\rm\cite{bcc08, c89,m92}}
Let $T$ be a tree, $(a_1,\ldots, a_n)$ is the size  of the maximum path partition of $\vec{T}_v$, then
$$f(T,v)=\sum_{i=1}^{n}2^{a_i}-n+1.$$
\end{thm}

\begin{lem}{\rm{(\cite{psv95})}}
The pebbling numbers of the cycles $C_{2n+1}$ and $C_{2n}$ are
$$f(C_{2n+1})=2\left\lfloor\frac{2^{n+1}}{3}\right\rfloor+1,f(C_{2n})=2^n.$$
\end{lem}

The following lemma is important in the proof of our main result.

\begin{lem}{\rm(\cite{m92}, No-Cycle-Lemma)}
If we have a graph $G$ with a certain distribution of pebbles
and a vertex $v$ of $G$ such that $m$ pebbles can be moved to $v$, then there always
exists an acyclic orientation $H$ for $G$ such that $m$ pebbles can still be moved
to $v$ in $H$.
\end{lem}

In this paper, we will show the pebbling number of Jahangir graphs $J_{n,m}$ when $n$ is even and $m\geq8$, which generalized the result of A. Lourdusamy etc.\cite{ljm11}.

\section{Main Result}

\begin{Def}
Jahangir Graph~$J_{n,m}$ {\rm($m\geq3$)} has $nm+1$ vertices, that is, a graph consisting of a cycle~$C_{nm}$ with one additional vertex which is
adjacent to $m$ vertices of $C_{nm}$ at distance $n$ to each other on $C_{nm}$.
\end{Def}

\begin{figure}
\begin{center}
\includegraphics[scale=1]{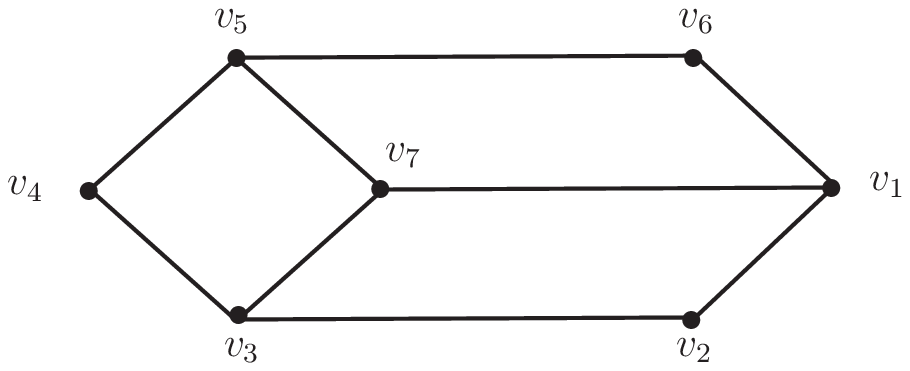}

\caption{{\small $J_{2,3}$.}\label{J23}}
\end{center}
\end{figure}

In this paper, we use the following notations: For a graph G, assume $v\in V(G)$, $A$ is a subset of $V(G)$, we use $p(v)$ and $p(A)$ to denote the number of pebbles on $v$ and $A$, respectively. We use $p_i$ to denote $p(v_i)$ for short.

Let $w$ be the root vertex of $G$. We say that a pebbling step from $u$ to $v$ is \emph{greedy} if $dist(v,w)<dist(u,w)$,
and that a graph $G$ is \emph{greedy} if from any distribution with $f(G)$ pebbles on $G$, one pebble can be moved to any specified root vertex $w$ with \emph{greedy} pebbling moves.

Clarke etc.\cite{chh97} asked the following question.

\begin{pps}{\rm\cite{chh97}}
Is every bipartite graph greedy?
\end{pps}

Here we give a counterexample.

\begin{thm}\label{thm2.3.0}
$J_{2,3}$ is not greedy.
\end{thm}

\begin{pf}
Assume $J_{2,3}$ is shown as Figure~{\rm\ref{J23}}, We know that $f(J_{2,3})=8$, if $p_2=p_6=3,p_1=p_7=1$, then we can not move a pebble on~$v_4$~with greedy pebbling moves.
\end{pf}

By the following construction, we can show a set of bipartite graphs that are not greedy.

\begin{lem}{\rm\cite{cpxy15,x15}}\label{lemcc}
Given a graph $G$ and a vertex $v$ on $V(G)$, form $G'$ from $G$ by adding a
new vertex, $u$, with $N(u)=N(v)$. If $G$ is Class 0, then $G'$ is also Class 0.
\end{lem}

\begin{cor}
Assume $J_{2,3}$ is shown as Figure~{\rm\ref{J23}}, Let $G_m$ be obtained from $J_{2,3}$ by adding $m$ new vertices, each of which is connected with all of $\{v_1,v_3,v_5\}$. Then $G_m$ is not greedy.
\end{cor}

\begin{pf}
By Lemma~{\rm\ref{lemcc}}, $G_m$ is Class 0, Let $p_2=p_6=3,p_3=p_4=p_5=0$, $p(v)=1$ otherwise, then we can not move a pebble on~$v_4$~with greedy pebbling moves. So $G_m$ is not greedy.
\end{pf}

A. Lourdusamy etc.\cite{ljm11} give the pebbling number of~$J_{2,m}$.

\begin{thm}{\rm\cite{ljm11}}\label{lemljm}
$f(J_{2,m})=2m+10$, $m\geq8$.
\end{thm}

In this paper, we will show the pebbling number of $J_{n,m}$ ($n$ is even, $m\geq8$).

Assume $J_{n,m}$ is obtained from $C_{nm}=v_0v_1\cdots v_{nm-1}$ with an additional vertex $u$ which is adjacent to~$v_{ni}$ $(0\leq i\leq m-1)$.
Let $P_i=v_{ni}v_{ni+1}\cdots v_{n(i+1)}$ for $i=0,1,\ldots,m-1$.

\begin{lem}\label{lem1}
Let $L_n=v_0v_1\cdots v_n$ be a path with length $n$, $D$ is a pebbling distribution on $L_n$, then

1) If neither of $v_0$ and $v_n$~are reachable, then $|D|\leq f(C_{n})-1$.

2) If only one of $v_0$ and $v_n$ is reachable, but not $2$-reachable, then $|D|\leq f(C_{n+1})-1$.

3) If both of $v_0$ and $v_n$ are reachable, but not $2$-reachable, then $|D|\leq f(C_{n+2})-1$.
\end{lem}

\begin{pf}
1)If neither of $v_0$ and $v_n$ are reachable, then we consider the same pebbling distribution $D$ on a new graph $C_{n-1}=v_0v_1\cdots v_{n-1}v_0$,
 So $D$ is not $v_0$-solvable on $C_{n}$, thus $|D|\leq f(C_{n})-1$.

2)If only one of $v_0$ and $v_n$ is reachable, but not $2$-reachable, without loss of generality, assume $D$ is $v_n$-solvable, then we joint $v_0$ and $v_n$ to get a new graph $C_{n+1}=v_0v_1\cdots v_{n}v_0$, $D$ is not $v_0$-solvable, so $|D|\leq f(C_{n+1})-1$.

3)If both of $v_0$ and $v_n$ are reachable, but not $2$-reachable, then we add a new vertex $u$ and new edges $uv_0$,$uv_n$ on $L_n$, to get a new graph $C_{n+2}=uv_0v_1\cdots v_nu$. $D$ is not $u$-solvable on $C_{n+2}$, thus $|D|\leq f(C_{n+2})-1$.
\end{pf}

By direct calculation, we can get

\begin{lem}\label{lem2}
$f(C_{n-1})+f(C_{n+1})\geq 2f(C_{n})$ for $n\geq2$.
\end{lem}

\begin{lem}\label{lem3}
Assume $D$ is a pebbling distribution on $J_{n,m}$ ($m\geq8$), if none of $u$, $v_0$ and $v_n$ are reachable, and $D(P_0)=0$, then we can get the following tight upper bounds.
\begin{align*}
|D|\leq\alpha=\left\{
\begin{array}{ll}
\frac{m}{2}(f(C_n)+f(C_{n+2})-2)-f(C_{n+2})+1, & \mbox{if}~m~\mbox{is even},\\
\frac{m-1}{2}(f(C_n)+f(C_{n+2})-2)+f(C_{n+1})-f(C_{n+2}), & \mbox{if}~m~\mbox{is odd}.
\end{array}
\right.
\end{align*}
\end{lem}

\begin{pf}
Note that $J_{n,m}$ divide the cycle $C_{nm}$ to paths $P_i=v_{ni}v_{ni+1}\cdots v_{n(i+1)}$ for $i=0,1,\ldots,m-1 $ with length $n$.

We only need to consider the contribution of the pebbling distribution on every path $P_i$ to its endpoints.

If the pebbles on $P_i$ can make both of its endpoints reachable, the number of pebbles is denoted by $L$ (Large);
If the pebbles on $P_i$ can  only make one of its endpoints reachable, the number of pebbles is denoted by $M$~(Middle);
If the pebbles on $P_i$ can make neither of its endpoints reachable, the number of pebbles is denoted by $S$~(Small). We also use $L,M$ or $S$ to name the path with $L,M$ or $S$ pebbles on it, respectively.

We will give a pebbling distribution with as many as possible pebbles.

\textbf{Claim:}
\begin{enumerate}
   \item Any two paths with~$L$~pebbles cannot be adjacent, $D(P_0)=0$, neither of $P_1$ and $P_{m-1}$ is $L$.
  \item At least one $S$ is between any two $L$.
Otherwise, there is a sequence $L,M,M,\ldots,M,L$, then one endpoint is reachable from two paths respectively.
  \item  At least one $S$ is between $L$ and $P_0$.
Otherwise, there is a sequence $0,M,M,\ldots,M,L$, then one endpoint is reachable from two paths  or the endpoint of $P_0$ is reachable from $L$ which is adjacent to $P_0$.
  \item We may assume that there does not exist two $M$ adjacent.
By Lemma~\ref{lem2}, if we replaced them by $S,L$, we can get a new distribution with at least the same number of pebbles.
\end{enumerate}

Assume there are $k$ paths with $L$ pebbles, from the Claim above, there are at least $k+1$ paths with $S$ pebbles, the number of pebbles on each of the left $m-2k-2$ paths is at most $M$.Thus~$|D|\leq kL+(k+1)S+(m-2k-2)M.$

If $m$ is even, then at most $\frac{m}{2}-1$ paths are $L$. By Lemma~\ref{lem2}, $|D|\leq kL+(k+1)S+(m-2k-2)M\leq kL+(k+1)S+\frac{m-2k-2}{2}(L+S)$.

If $m$ is odd, then at most $\frac{m-3}{2}$ paths are $L$. By Lemma~\ref{lem2}, $|D|\leq kL+(k+1)S+(m-2k-2)M\leq kL+(k+1)S+\frac{m-2k-3}{2}(L+S)+M$.

The upper bounds of $L$,$M$,$S$ are given by Lemma~\ref{lem1}, and we are done.

Now we give the distribution $D^*$ to show that the bounds are tight:

For the path sequence $P_0,P_1,\ldots,P_{m-1}$,

If $m$ is even, the sequence of distribution is $0,S,L,S,L,\ldots,S$;

If $m$ is odd, the sequence of distribution is
$0,S,L,\ldots,S,L,M,S$, where $M$ is $v_{n(m-1)}$-solvable.

Let $L$,$M$ and $S$ get the upper bounds in Lemma~\ref{lem1}, we are done.
\end{pf}\vskip .5cm

Let $\varsigma$ be a sequence of pebbling moves from the distribution $D_1$ to $D_2$ on the graph $G$, then we say that \emph{the number of the pebbles cost} in $\varsigma$ is $|D_1|-|D_2|$.
The pebbling move along one edge $\{w,v\}$ from $w$ to $v$ induce a directed edge $(w,v)$, similarly, we can define the directed graph $\vec{G}$ induced by $\varsigma$, in which we allow some edges have no direction. \emph{The source vertex} of $\vec{G}$ is the vertex that its out-degree is greater than 0, and its in-degree is 0;  \emph{The sink vertex} of $\vec{G}$ is the vertex that its in-degree is greater than 0, and its out-degree is 0 (the vertex with no directed edges has the out-degree 0 and in-degree 0).

Assume $w$ and $v$ are adjacent, the sequence of pebbling moves $\eta$ remove $2\beta$ pebbles from $w$, and place $\beta$ pebbles on $v$, then we can say that \emph{$\eta$ moves $\beta$ pebbles from $w$ to $v$}. Now we generalize it to two nonadjacent vertices. We paint the pebbles on $w$ red, paint the pebbles on other vertices black. If a pebbling move remove two pebbles from $v_1$, and place one pebble on $v_2$. We consider three cases before this pebbling move.

(1) There are at least two red pebbles on $v_1$, then we remove two red pebbles from $v_1$, and place one red pebble on $v_2$;

(2) There is only one red pebble on $v_1$, then we remove one red and one black pebble from $v_1$, and place one red pebble on $v_2$;

(3) There is no red pebble on $v_1$, then we remove two black pebbles from $v_1$, and place one black pebble on $v_2$.

If $\gamma$ pebbles on $v$ are red after a sequence of pebbling moves $\varphi$, then we say that \emph{$\varphi$ moves $\gamma$ pebbles from $w$ to $v$}.

\begin{thm}\label{thm1}
Let $n$ be an even integer, $m\geq8$,
$f(J_{n,m})=f_{2^{\frac{n}{2}+1}-1}(C_{n+2})+\alpha+1$, where $\alpha$ is shown in Lemma~\ref{lem3}.
\end{thm}

\begin{pf}
First we show that $f(J_{n,m},v_{\frac{n}{2}})=f_{2^{\frac{n}{2}+1}-1}(C_{n+2})+\alpha+1$, assume the target vertex is $v_{\frac{n}{2}}$.

\textbf{Lower bound}:
Let the distribution $D$ be obtained from the distribution $D^*$ given in the proof of Lemma~\ref{lem3}  by adding $f_{2^{\frac{n}{2}+1}-1}(C_{n+2})$ pebbles on $v_{4n+n/2}$.

We only need to show that $D$ is not $v_{\frac{n}{2}}$-solvable. We show it by contradiction. If $D$ was $v_{\frac{n}{2}}$-solvable, then there exist a sequence of pebbling moves so that one pebble can be moved from $D$ to $v_{\frac{n}{2}}$. Assume $\tau$ is such a sequence of pebbling moves with the number of the pebbles cost minimum.

According to the No-Cycle-Lemma, the directed graph $\vec{J}$ induced by $\tau$ has no directed cycle. For $\tau$ is the sequence of pebbling moves that cost the minimum number of pebbles, $\vec{J}$ has only one sink vertex $v_{\frac{n}{2}}$. If the out-degree of some vertex $v$ is not 0, then $\tau$ moves one pebble from $v$ to $v_{\frac{n}{2}}$.

Moreover, we claim the following:
\begin{itemize}
  \item  $(u,v_{ni})\notin\vec{J}$ for $2\leq i\leq m-1$.

   Assume that $(u,v_{nj})\in\vec{J}$ for some $2\leq j\leq m-1$, we may assume that one pebble have been moved from $u$ to $v_{nj}$. We paint this pebble red, and we paint other pebbles black. note that $\vec{J}$ has only one sink vertex $v_{\frac{n}{2}}$, so the red pebble must be moved along $P_{j-1}$ to $v_{n(j-1)}$ (or along $P_{j+1}$ to $v_{n(j+1)}$). Then we choose the pebbling move from $u$ to $v_{n(j-1)}$ (or $v_{n(j+1)}$) instead of this sequence of pebbling moves to get a new sequence of pebbling moves with less pebbles cost than $\tau$, which is a contradiction to the assumption that $\tau$ costs the least pebbles.

  \item $(v_{4n},v_{4n-1})\notin \vec{J}$ and $(v_{5n},v_{5n+1})\notin\vec{J}$.

  Note that for the path sequence $P_0,P_1,P_2,P_3,P_4,P_5$, the sequence of pebbles is$0,S,L,S,$\newline$L',S$, where $L'$ is obtained from $L$ by adding $f_{2^{\frac{n}{2}+1}-1}(C_{n+2})$ pebbles on $v_{4n+n/2}$. We paint the pebbles on $P_4$ red, and paint the other pebbles black. Assume $(v_{4n},v_{4n-1})\in \vec{J}$, then one red pebble has been moved from $v_{4n}$ to $v_{4n-1}$, according to the No-Cycle-Lemma, and the assumption that  $\tau$ costs the least pebbles, we can get that this red pebble must be moved along $P_3$ to $v_{3n}$, now we can get $(v_{3n},u)\notin \vec{J}$, for if $(v_{3n},u)\in \vec{J}$, then the red pebble is moved to $u$, so we choose the pebbling move from $v_{4n}$ to $u$ instead of such sequence of pebbling moves, and we get a new sequence of pebbling moves with  less pebbles cost than $\tau$, which is a contradiction. From Claim 1), we can get $(u,v_{3n})\notin \vec{J}$, We know that $v_{3n}$ is not a sink vertex in $\vec{J}$, so we can get $(v_{3n},v_{3n-1})\in\vec{J}$. By a similar argument, we can get that the red pebble must be moved along the path sequence $P_3,P_2,P_1$ to $v_n$. It cost us at least five red pebbles to move one red pebble to $v_n$, so we can choose a new sequence of pebbling moves instead of this sequence: remove four red pebbles from $v_{4n}$, and add two red pebbles on $u$, and then remove these two red pebbles and add one red pebble to $v_n$, which is a contradiction. Similarly, we can show that $(v_{5n},v_{5n+1})\notin\vec{J}$.
\end{itemize}

 According to the claim above, we can get that the pebbles must be moved from $P_4$ to $u$, then from $u$ to $P_0$, directly. At most $2^{\frac{n}{2}+1}-1$ red pebbles can be moved to $u$, so we can not move one red pebble to $v_{\frac{n}{2}}$.

\textbf{Upper bound}:
Assume~$f_{2^{\frac{n}{2}+1}-1}(C_{n+2})+\alpha+1$ pebbles are placed on $J_{n,m}$. Let $C^i=P_i\bigcup u$ be the cycle induced by $P_i$ and $u$ in $J_{n,m}$.

We first consider the case: $p(P_0)=0$.

We only need to show that with $f_{t-1}(C_{n+2})+\alpha+1$ pebbles on $J_{n,m}$ such that $p(P_0)=0$, one can move $t$ pebbles to $C^0$ without the pebbles on $S$. So we may assume that $p(u)=0$. We use induction on $t$. It holds for $t=1$ by Lemma~\ref{lem3}. Assume it holds for $t-1$ ($t\geq 2$). For the case $t$, by Lemma~\ref{lem3}, one of the following holds:
1)Two $L$ are adjacent; 2) one $L$ and one $M$ are adjacent, so that two pebbles can be moved to the joint vertex, and so one pebble can be moved to $u$;
3) The number of the pebbles on some path $P_j$ is larger than $f(C_{n+2})$. 

Case 1. Two $L$ are adjacent. Assume that $P_k$ and $P_{k+1}$ both have large number of pebbles($L$). We can move one pebble to $v_{nk}$ from each path, and so one pebble can be moved to $u$.
Then we replace the distribution of pebbles left on $P_k$ and $P_{k+1}$ by the distribution with $2^{n/2}-1$ pebbles on $P_k$ and $P_{k+1}$ respectively, such that none of its endpoints are reachable (just put $2^{n/2}-1$ pebbles on $v_{nk+n/2}$ and $v_{n(k+1)+n/2}$, respectively). Then the total number of the new distribution of $J_{n,m}\backslash u$ is
\begin{align*}
&f_{t-1}(C_{n+2})+\alpha+1-|P_k|-|P_{k+1}|+2(2^{n/2}-1)\\
\geq& f_{t-1}(C_{n+2})+\alpha+1-2(2^{n/2+1}-1)+2(2^{n/2}-1)\\
=&f_{t-2}(C_{n+2})+\alpha+1.
\end{align*}
By induction, $t-1$ pebbles can be moved to $u$ from the new distribution, and we do not use the pebbles on $P_k$ or $P_{k+1}$, that means we can move $t-1$ pebbles from the original distribution, and we are done.

Case 2. The proof is similar to Case 1.

Case 3. The number of the pebbles on some path $P_j$ is larger than $f(C_{n+2})$, we can move one pebble to $u$ from $P_j$ with at most $2^{n/2+1}-1$ pebbles cost (for the even cycle is greedy). Then we left at least $f_{t-2}(C_{n+2})+\alpha+1$ pebbles on $J_{n,m}\backslash u$, and we can move $t-1$ pebbles to $u$ with the left pebbles by induction, and we are done.

Now if $p(P_0)>0$, then we only need to move $2^{n/2+1}-p(P_0)$ pebbles from $J_{n,m}\backslash \{P_0\bigcup u\}$ to $u$ (then we have $2^{n/2+1}$ pebbles on $P_0\bigcup u$, so one pebble can be moved to $v_{n/2}$). The number of pebbles on $J_{n,m}\backslash \{P_0\bigcup u\}$ is  $f_{2^{\frac{n}{2}+1}-1}(C_{n+2})+\alpha+1-p(P_0)$, which is larger than $f_{2^{n/2+1}-p(P_0)-1}(C_{n+2})+\alpha+1$. So we are done from the argument above.

If the target vertex is not $v_{n/2}$, we only need to check the upper bound, we may assume that the target vertex belongs to $P_0$, then by a similar argument, one can show that $2^{n/2+1}$ pebbles can be moved to $P_0\bigcup u$, and so
one pebble can be moved to the target vertex with $f_{2^{\frac{n}{2}+1}-1}(C_{n+2})+\alpha+1$ pebbles on $J_{n,m}$, and we are done.
\end{pf}

\section{Remark}

Let $n=2$, then we can get $f(J_{2,m})=2m+10$ for $m\geq 8$ from Theorem~\ref{thm1}, which is just Lemma~\ref{lemljm}. For $m<8$, there may exist a sequence of pebbling moves with least pebbles cost, which does not through the vertex $u$, so we cannot get the tight lower bound with the method in Theorem~\ref{thm1}, but we can still get the upper bound by a similar argument.

\end{document}